\documentclass  [journal,onecolumn,tbtags, draftcls, 12pt]{IEEEtran}

\usepackage[sort,compress]{cite}
\usepackage{graphicx}

\usepackage[cmex10]{mathtools}
\usepackage{dsfont,amssymb,bm}
\usepackage{arydshln}

\usepackage[standard,amsmath,thmmarks] {ntheorem}
\newtheorem* {problem}       {Problem (P1)}

\usepackage[T1]{fontenc}

\let\VEC      \relax
\DeclareMathOperator\VVEC{vec}
\DeclareMathOperator\DIAG{diag}
\let\MAT      \mathbf

\def\tMAT#1{\tilde{\MAT #1}}

\DeclareMathOperator*{\TR}{Tr}

\newcommand\PR {\mathds{P}}
\newcommand\EXP{\mathds{E}}

\newcommand\reals{\mathds{R}}

\newcommand\important[1]{\textbf{\textit{ #1}}}

\newcommand\NORMAL{\mathcal{N}}
\newcommand\TRANS{\intercal}
\newcommand\COM{\mathit{com}}

\begin{document}

\title {Sufficient statistics for linear control strategies in
decentralized systems with partial history sharing}

\author{Aditya Mahajan and Ashutosh Nayyar}

\maketitle

 \begin{abstract}
   In decentralized control systems with linear dynamics, quadratic cost, and
   Gaussian disturbance (also called decentralized LQG systems) linear control strategies
   are not always optimal. Nonetheless, linear control strategies are appealing
   due to analytic and implementation simplicity. In this paper, we investigate
   decentralized LQG systems with partial history sharing information structure and
   identify finite dimensional sufficient statistics for such systems. Unlike
   prior work on decentralized LQG systems, we do not assume partially
   nestedness or quadratic invariance. Our approach is based on the common
   information approach of Nayyar \emph{et al}, 2013 and exploits the linearity
   of the system dynamics and control strategies. To illustrate our methodology, we identify sufficient statistics for linear strategies in decentralized systems where controllers communicate over a strongly connected graph with finite delays, and for decentralized systems consisting of coupled subsystems with control sharing or one-sided one step delay sharing information structures.
 \end{abstract}

\section{Introduction}

With the increasing applications of networked control systems, the problem of
finding the best linear control strategy for decentralized systems with linear
dynamics, quadratic cost, and Gaussian disturbances (henceforth referred to as
decentralized LQG systems) has received considerable attention in recent
years~\cite{MMRY:tutorial-CDC} (and references therein).

In centralized LQG systems, linear control strategies are globally
optimal, the best linear control strategies are characterized by the solution of a
Riccati equation, the best linear control is a function of the
controller's estimate of the state of the plant and 
this estimate is updated using Kalman filtering equations.
In contrast, the problem of finding the best linear control strategies for
decentralized LQG systems has the following salient features:
\begin{enumerate}
  \item In general, linear control strategies are not globally optimal, i.e.,
    there may exist non-linear control strategies that outperform linear
    strategies as is illustrated by the Witsenhausen
    counterexample~\cite{Witsenhausen:1968} and memoryless control in Gaussian
    noise~\cite{LipsaMartins:2011b}. Linear strategies are globally optimal only
    when the controller has specific information structure such as
    static~\cite{Radner:1962}, partially nested~\cite{HoChu:1972}, or
    stochastically nested~\cite{Yuksel:2009} information structures and their
    variations.

  \item In general, the problem of finding the best linear control strategies is
    not convex. It may be converted to a convex model matching problem only 
    when the sparsity pattern of the plant and the controller
    have specific structure such as funnel causality~\cite{BamiehVoulgaris:2005}
    or quadratic invariance~\cite{RotkowitzLall:2006} and their variations.

  \item In general, the best linear control strategy may not have a finite
    dimensional sufficient statistic, i.e., it may not be possible to represent
    the best linear controller by a finite set of estimates that are generated
    by recursions of finite order as is illustrated by the two controller
    completely decentralized system considered in~\cite{WhittleRudge:1974}. The
    best linear strategies are known to have a finite dimensional sufficient
    statistic only for specific examples~\cite{Yoshikawa:1975,
      VaraiyaWalrand:1978, SwigartLall:2010, ShahParrilo:2013, LessardLall:2011,
    KhorsandAlamGattami:2012, NayyarKalathilJain:2013, LessardNayyar:2013}. Note
    that all of these examples have partially nested information structure and
    some of these examples have quadratic invariant sparsity structure. It is
    generally believed that the best linear control strategies in partially
    nested and quadratic invariant systems will have finite dimensional
    sufficient statistic.
\end{enumerate}

In this paper, we investigate the third aspect of decentralized LQG systems
described above, viz., whether finite dimensional sufficient statistics for linear control strategies can be identified for some subclass of decentralized LQG systems. In
particular, we investigate decentralized LQG systems with partial history
sharing information structure~\cite{NMT:partial-history-sharing}, which is a
generalization of several well-known information structures of decentralized
control. The partial history sharing model, in general, is not partially nested
or quadratic invariant. Our main results for this model are presented in
Section~\ref{sec:results} and and can be summarized as follows:
\begin{enumerate}
  \item we identify finite dimensional sufficient statistics for the best linear control
    strategy; and
  \item we show that the update equation of these sufficient statistics is
    similar to Kalman filter updates.
\end{enumerate}

In Section~\ref{sec:examples}, we apply these results to decentralized control
systems in which the controllers communicate along a strongly connected graph
with finite delay between any pair of controllers. In Section~\ref{sec:others},
we show that these results can be also applied to models that are not partial
history sharing, but can be converted to one by using a person-by-person
approach.

To the best of our knowledge, \important{these are the first sufficient statistics results
for best linear strategies in decentralized LQG systems that are neither
partially nested nor quadratic invariant.} Our results suggest that the form of
the sufficient statistics is a consequence of linearity of system dynamics and control strategies
rather than partial nestedness or quadratic invariance of the information structure.

Our solution methodology is based on the \emph{common information
approach} developed in \cite{Nayyar:PhD} and used in~\cite{NMT:partial-history-sharing} for decentralized control systems with partial history sharing. However, our results cannot be derived
directly using the results of~\cite{NMT:partial-history-sharing}. For a general decentralized system with partial history sharing, the results
of~\cite{NMT:partial-history-sharing} provide the structure of globally optimal
control strategies \emph{and} a dynamic programming decomposition. In this paper, we exploit linearity (of control strategies and of the underlying decentralized system) to address only the problem of finding the structure of best linear strategy. We do not address the problem of computing the best linear strategy. This narrower focus allows us to get simpler results than in \cite{NMT:partial-history-sharing}.

Even with finite dimensional sufficient statistics, the problem of computing the best linear strategies is, in general, a non-convex optimization problem unless the system is quadratically 
invariant;  and even if the best linear strategy is identified, it is
globally optimal only if the system is partially nested. Nonetheless, when the system is either partially nested or quadratic invariant,
it may be possible to use finite dimensional sufficient statistics to compute best linear or globally optimal strategies. For example, an approach similar to ours was used
in~\cite{LessardNayyar:2013} to identify sufficient statistics for best linear
control strategies (that were also globally optimal) for a two player
decentralized LQG team~\cite{LessardNayyar:2013} that is partially nested and
quadratic invariant. The authors of~\cite{LessardNayyar:2013} then  exploited the
partial nested nature of the system to identify explicit expressions for the best linear control
strategies. 

\subsection*{Notation}

Uppercase letters denote random variables/vectors and lowercase letters denote
their realization. Bold uppercase letters denote matrices. $\PR(\cdot)$ denotes
the probability of an event and $\EXP[\cdot]$ denotes the expectation of a
random variable. $\reals$ denotes the set of real numbers. 

For a sequence of (column) vectors $X$, $Y$, $Z$, \dots, the notation $\VVEC(X,Y,Z,\dots)$
denotes the vector $[X^\TRANS, Y^\TRANS, Z^\TRANS]^\TRANS$. The vector
$\VVEC(X_1, \dots, X_t)$ is also denoted by $X_{1:t}$. 

The notation $\MAT A = \DIAG(\MAT B, \MAT C, \MAT D)$ denotes a block diagonal
matrices with blocks $\MAT B$, $\MAT C$, and $\MAT D$ on the diagonal. $\MAT
A^\TRANS$ denotes the transpose of a matrix and $\TR[\MAT A]$ denotes the trace
of a matrix. 

The notation $\MAT 0_{n \times m}$ denotes a $n \times m$ all zeros matrix;
$\MAT I_n$ denotes a $n \times n$ identity matrix. We omit the subscripts when
dimensions can be inferred from context.

For any two random vectors $X$ and $Y$, we say that $X$ is a sub-vector of $Y$,
and denote it by $X \subset Y$ if the set of all components of $X$ is a subset
of the set of all components of $Y$. More formally, $X \subset Y$ if there
exists a row-stochastic binary matrix $\MAT P$ (i.e., all its elements are $0$
or $1$ and each row has a single $1$) such that $X = \MAT P Y$. 

\section{Problem formulation}

\subsection{Model} \label{sec:model}

Consider a linear dynamic system with $n$ controllers and a partial history sharing
information structure~\cite{NMT:partial-history-sharing}. We follow the same notation
as~\cite{NMT:partial-history-sharing} and, for completeness, restate the model
below.

The system operates in discrete
time for a horizon $T$. Let $X_t \in \reals^{d_x}$ denote the state of the
system at time $t$, $U^i_t \in \reals^{d^i_u}$ denote the control action of
controller~$i$, $i=1,\dots,n$ at time $t$, and $\VEC U_t$ denote the vector
$\VVEC(U^1_t, \dots, U^n_t)$.

The initial state $X_1$ has a probability distribution $\NORMAL(0,\Sigma_x)$ and
evolves according~to
\begin{equation}
  \label{eq:state}
    X_{t+1} = \MAT A_t X_t + \MAT B_t \VEC U_t + W^0_t
\end{equation}
where $\MAT A_t$ and $\MAT B_t$ are matrices of appropriate dimensions and
$\{W^0_t\}_{t=1}^{T}$ is a sequence of i.i.d\@. zero-mean Gaussian random
variables with probability distribution~$\NORMAL(0,\Sigma_{w^0})$.

As in~\cite{NMT:partial-history-sharing}, at any time $t$, each controller has
access to three types of data: the current observation~$Y^i_t$, the local
memory~$M^i_t$, and the shared memory~$\VEC Z_{1:t-1}$. The details of the
information structure will be described later. We use $Y_t$ to denote $\VVEC(Y^1_t,\ldots, Y^n_t)$ and  $M_t$ to denote $ \VVEC(M^1_t,\ldots, M^n_t)$.

We restrict attention to linear control strategies and assume that
controller~$i$'s strategy is of the form:
\begin{equation}
  \label{eq:linear-control}
  U^i_t = \MAT K^i_t \VEC Z_{1:t-1} + \MAT G^i_t Y^i_t + \MAT H^i_t M^i_t
\end{equation}
where $\MAT K^i_t$, $\MAT G^i_t$, and $\MAT H^i_t$ are matrices of appropriate
dimensions. The collection of $\{(\MAT K^i_t, \MAT G^i_t, \MAT H^i_t)\}_{t=1}^T$
is referred to as the  \emph{control strategy} of controller~$i$. 

Combining \eqref{eq:linear-control} for all controllers, we can write
\begin{equation}
  \label{eq:linear-control2}
  U_t = \MAT K_t \VEC Z_{1:t-1} + \MAT G_t Y_t + \MAT H_t M_t,
\end{equation}
where  $\MAT K_t = [\MAT K^{1\,\TRANS}_t \mid \dots \mid \MAT
K^{n\,\TRANS}_t]^\TRANS_{\strut}$, $\MAT G_t = \DIAG(\MAT G^1_t, \dots, \MAT G^n_t)$ and $\MAT
H_t = \DIAG(\MAT H^1_t, \dots, \MAT H^n_t)$. 

At time~$t$, the system incurs a quadratic cost $\ell(X_t, \VEC U_t)$ given by
\begin{equation}
  \label{eq:cost}
  \ell(X_t, \VEC U_t) = X_t^\TRANS \MAT Q X_t + \VEC U_t^\TRANS \MAT R \VEC U_t
\end{equation}
where $\MAT Q$ is positive semi-definite and $\MAT R$ is positive definite
matrices of appropriate dimensions.

We are interested in choosing control strategies of all controllers to minimize
\begin{equation}
  \label{eq:total-cost}
  \EXP
  \Big[ \sum_{t=1}^T \ell(X_t, \VEC U_t) \Big],
\end{equation}
where the expectation is with respect to the joint probability measure on
$(X_{1:T}, \VEC U_{1:T})$ induced by the choice of the control strategies.

\subsection{Partial history sharing information structure}

As described earlier, controller $i$ has access to three types of data at
time~$t$: the current observation~$Y^i_t$, the local memory~$M^i_t$, and the
shared memory~$\VEC Z_{1:t-1}$. These variables are given as follows:
\begin{enumerate}
  \item The current local observation $Y^i_t \in \reals^{d^i_y}$ of
    controller~$i$ is given by
    \begin{equation}  \label{eq:observation}
      Y^i_t = \MAT C^i_t X_t + W^i_t
    \end{equation}
    where $\MAT C^i_t$ is a matrix of appropriate dimensions and $\{W^i_t\}_{t=1}^{T}$
    is a sequence of i.i.d\@. zero-mean Gaussian random variables with probability distribution
    $\NORMAL(0,\Sigma_{w^i})$. The random variables in the collection
    $\{X_1,W^j_t,t=1,\dots,T, j =0, 1,\dots, n\}$, called \emph{primitive random
    variables},  are mutually  independent. Combining \eqref{eq:observation} for all controllers, we can write 
    \[ Y_t = \MAT C_t X_t + W^{1:n}_t,
    \]
where $W^{1:n}_t = \VVEC(W^1_t,\ldots,W^n_t)$ and $\MAT C_t = [\MAT C^{1\,\TRANS}_t \mid \dots \mid \MAT
C^{n\,\TRANS}_t]^\TRANS_{\strut}$.
  \item The local memory $M^i_t \in \reals^{d^i_m}$ of controller~$i$ is a
    subvector of the history of its local observations and actions:
    \begin{equation}\label{eq:memory}
      M^i_t \subset \{Y^i_{1:t-1},U^i_{1:t-1}\}
    \end{equation}
     At $t=1$, the local memory is empty, which we will represent by the convention $M^i_1 := 0$.

  \item In addition, all controllers have access to a shared memory $\VEC Z_{1:t-1}$, where $Z_t
    = \VVEC(Z^1_t, \dots, Z^n_t)$. The shared memory $\VEC Z_{1:t-1}$ 
    is a subset of the history of observations and actions of all
    controllers:
    \begin{equation}\label{eq:shared memory}
      \VEC Z_{1:t-1} \subset \{\VEC Y_{1:t-1}, \VEC U_{1:t-1}\}.
    \end{equation}
At $t=1$, the shared memory is
    empty, $\VEC Z_0 :=0$; at each time $Z_t \in \reals^{d_z}$.
\end{enumerate}

The local and shared memories are updated as follows: After taking the control action at time~$t$, controller~$i$ sends a subvector
    $Z^i_t$ of its local information $\{M^i_t, Y^i_t, U^i_t\}$ to the shared memory. We assume that the protocol of choosing the subset $Z^i_t$ is pre-specified.
 After sending data~$Z^i_t$ to the shared memory, controller~$i$ updates its local
memory according to a pre-specified protocol such that $M^i_{t+1} \subset
\{M^i_t,Y^i_t,U^i_t\}\setminus Z^i_t$, which ensures that the contents of the
local and shared memories do not overlap. 

The process of generating   the new local memory $M^i_{t+1}$ and $Z^i_t$ described above can be written in terms of the following equations:
\begin{equation}
  \label{eq:update_1}
  M^i_{t+1}
  = \MAT P^i_{mm,t}M^i_t +  \MAT P^i_{my,t}Y^i_t + \MAT P^i_{mu,t} U^i_t    
\end{equation}
and
\begin{equation}
  \label{eq:update_2}
  Z^i_{t}
   = \MAT P^i_{zm,t}M^i_t +  \MAT P^i_{zy,t}Y^i_t + \MAT P^i_{zu,t} U^i_t , 
\end{equation}
where $\MAT P^i_{**,t}$ are  matrices  that satisfy the following properties:
\begin{enumerate}
\item[A1.] Each entry of $\MAT P^i_{**,t}$ is either $0$ or $1$.
\item[A2.] The matrix 
\begin{equation*}
\begin{bmatrix}
      \MAT P^i_{mm,t} & \MAT P^i_{my,t} &\MAT P^i_{mu,t} \\
      \MAT P^i_{zm,t} & \MAT P^i_{zy,t} &\MAT P^i_{zu,t}
\end{bmatrix}
\end{equation*}
is doubly stochastic (that is, each row and column sum is $1$).
\end{enumerate}

Note that the $\MAT P^i_{**,t}$ matrices are specified a priori based on the
memory update protocols of the system. Also note that properties A1 and A2 are a
consequence of these memory update protocols. We refer the reader to \cite{NMT:partial-history-sharing} for several examples of partial history sharing information structures.

Combining \eqref{eq:update_1} for all controllers we get
\begin{equation}
  \label{eq:update_3}
  M_{t+1}
    = \MAT P_{mm,t}M_t +  \MAT P_{my,t}Y_t + \MAT P_{mu,t} U_t    
\end{equation}
where $\MAT P_{mm,t} = \DIAG(\MAT P^1_{mm,t}, \dots, \MAT
    P^n_{mm,t})$,  $\MAT P_{my,t} = \DIAG(\MAT P^1_{my,t}, \dots, \MAT
    P^n_{my,t})$, $\MAT P_{mu,t} = \DIAG(\MAT P^1_{mu,t}, \dots, \MAT
    P^n_{mu,t})$.
    
     Similarly, combining \eqref{eq:update_2} for all controller gives
\begin{equation}
  \label{eq:update_4}
  Z_{t}
      = \MAT P_{zm,t}M_t +  \MAT P_{zy,t}Y_t + \MAT P_{zu,t} U_t    
\end{equation}
where $\MAT P_{zm,t} = \DIAG(\MAT P^1_{zm,t}, \dots, \MAT
    P^n_{zm,t})$,  $\MAT P_{zy,t} = \DIAG(\MAT P^1_{zy,t}, \dots, \MAT
    P^n_{zy,t})$, $\MAT P_{zu,t} = \DIAG(\MAT P^1_{zu,t}, \dots, \MAT
    P^n_{zu,t})$.
    
    An example of the above model is the delayed sharing information
structure~\cite{NMT:delay-sharing}, in which the shared memory consists of $k$
steps old observations and control actions all controllers, i.e., $\VEC
Z_{1:t-1} = \VVEC(Y_{1:t-k}, U_{1:t-k})$ and the local memory consists of the
observations and actions taken at $t-k+1$, \ldots, $t-1$, i.e., $M^i_t =
\VVEC(Y^i_{t-k+1:t-1}, U^i_{t-k+1:t-1})$. In particular, when the delay $k=2$,
then $M^i_t = \VVEC(Y^i_{t-1}, U^i_{t-1})$, $Z^i_t = \VVEC(Y^i_{t-2},
U^i_{t-2})$, and the equations for  generating   $M^i_{t+1}$ and $Z^i_t$ can be written as
\begin{equation*}
  M^i_{t+1}
  = \MAT 0 M^i_t + \begin{bmatrix} \MAT I \\ \MAT 0 \end{bmatrix}Y^i_t + \begin{bmatrix} \MAT 0 \\ \MAT I \end{bmatrix} U^i_t    
\end{equation*}
and
\begin{equation*}
  Z^i_{t}
   = \MAT I M^i_t +  \MAT 0 Y^i_t + \MAT 0  U^i_t .
\end{equation*}
%
%
%

\subsection{Generalized partial history sharing information structure}

We now describe the generalized version of the partial history sharing
information structure. As in the original partial history sharing model,
controller $i$ has access to three types of data at time~$t$: the current
observation~$Y^i_t$,  a shared memory $\VEC Z_{1:t-1}$ that is available to all
controllers and local memory~$M^i_t$ with  $M^i_1 :=0$ and $Z_0:=0$. The
difference between the original model and the generalized one lies in the memory
update rules. In the partial history sharing model, the local  and shared
memories are updated according to \eqref{eq:update_3} and \eqref{eq:update_4},
where $\MAT P_{**,t}$ are block diagonal matrices and $\MAT P^i_{**,t}$ satisfy
properties A1 and A2. In \emph{generalized partial history sharing information
structure}, the local  and shared memory update rules still satisfy
\eqref{eq:update_3} and \eqref{eq:update_4}, but we allow  $\MAT P_{**,t}$  to
be arbitrary matrices. We will describe examples of this information structure
in Section~\ref{sec:examples}.

\begin{remark} \label{rem:memory}
  In some cases, the
  local memory $M^i_t$ is always empty. In such systems, the update
  equations~\eqref{eq:update_1}-\eqref{eq:update_4} can be replaced by 
  \begin{equation} \label{eq:update-2}
    Z^i_{t} =  \MAT P^i_{zy,t}Y^i_t + \MAT P^i_{zu,t} U^i_t , 
  \end{equation}
  \begin{equation}
    Z_{t} =   \MAT P_{zy,t}Y_t + \MAT P_{zu,t} U_t    
    \end{equation}
\end{remark}
\subsection{Problem formulation}

We are interested in the problem of finding the best linear control strategies.
Specifically:

\begin{problem}
  For the model described above, given horizon~$T$, the matrices $\MAT A_t$,
  $\MAT B_t$, $\MAT C^i_t$, $\MAT Q$, $\MAT R$,
  the covariance matrices $\Sigma_{x}$, $\Sigma_{w^i}$, and the
  protocols for updating the local and shared memory, 
  find a control strategy of the form~\eqref{eq:linear-control}
  that minimizes the expected total  cost given by~\eqref{eq:total-cost}.
\end{problem}

One of the difficulties for Problem~(P1) is that the shared memory $Z_{1:t-1}$
available to all controllers is increasing with time; consequently, the size of the gain
matrices $\MAT K_t$ in \eqref{eq:linear-control2} is increasing as well . We identify appropriate sufficient
statistics $\breve X_t$ (to be defined later) that have the same dimension as
$\VVEC(X_t, Y_t, M_t)$ and show that the optimal controller is of the form
\[
    U_t = \tMAT K_t \breve X_t+ \MAT G_t Y_t + \MAT H_t M_t.
  \]
Furthermore, $\breve X_t$ may be updated in a manner similar to Kalman filtering
updates.

\section{Main results} \label{sec:results}

\subsection{A sub-problem and the induced centralized system}

The main idea of the proof is as follows. Arbitrarily fix the matrices $(\MAT
G_{1:T}, \MAT H_{1:T})$.
 Consider the sub-problem of finding the best
choice of matrices $\MAT K_{1:T}$  to minimize the total expected cost given by~\eqref{eq:total-cost}.

Following~\cite{NMT:partial-history-sharing}, we introduce a new decision
maker---the \emph{coordinator}---that sequentially observes the process
$\{Z_t\}_{t=1}^T$ and chooses actions $\tilde U_t = \VVEC(\tilde U^1_t, \dots,
\tilde U^n_t)$ where
\begin{equation}
  \label{eq:coordinator-control}
  \tilde U^i_t = \MAT K^i_t \VEC Z_{1:t-1}.
\end{equation}
The controllers of the original system are passive agents that generate~$U^i_t$
according to
\begin{equation}
  \label{eq:coordinator-passive}
  U^i_t = \tilde U^i_t + \MAT G^i_t Y^i_t + \MAT H^i_t M^i_t.
\end{equation}

Combine~\eqref{eq:coordinator-control} and~\eqref{eq:coordinator-passive} in
vector form to write
\begin{align}
  \tilde U_t &= \MAT K_t Z_{1:t-1}, 
  \label{eq:vector-control} \\
  U_t &= \tilde U_t + \MAT G_t Y_t + \MAT H_t M_t;
  \label{eq:vector-passive}
\end{align}
where $\MAT G_t$ and $\MAT H_t$ are block diagonal matrices and $\MAT
K_t$ is a stacked matrix as defined earlier.

As in~\cite{NMT:partial-history-sharing}, the optimization problem at the
coordinator is equivalent to a partially observed centralized stochastic control
problem, which we call the \emph{coordinated system}. Define the state $\tilde
X_t$ and the observation $\tilde Y_t$ of this coordinated system as:
\begin{align}
  \tilde X_t &= \VVEC(X_{t}, \VEC Y_{t}, \VEC M_{t}), \\
  \tilde Y_t &= \VEC Z_{t-1}.
\end{align}
Then the control action $\tilde U_t$ of this system is chosen according
to~\eqref{eq:vector-control} which is a linear functional of the observation
history. 

The coordinated system is a centralized system LQG system with linear dynamics, linear
observations, quadratic cost, and Gaussian disturbance. In particular:
\begin{enumerate}
  \item 
    The coordinated system has linear dynamics which may be written as
    \begin{equation} \label{eq:coord_dynamics}
      \tilde X_{t+1} = \tilde {\MAT A}_t \tilde X_t + \tilde {\MAT B}_t \tilde U_t +
      \tilde {\MAT F}_t \VEC W_t
    \end{equation}
    where $\VEC W_t = \VVEC(W^0_t, W^1_{t+1}, \dots, W^n_{t+1})$, and $\tilde {\MAT
    A}_t$, $\tilde {\MAT B}_t$, and $\tilde {\MAT F}_t$ are matrices of appropriate
    dimensions that are obtained by combining~\eqref{eq:state},
    \eqref{eq:observation}, \eqref{eq:update_3}, \eqref{eq:vector-control},
    and~\eqref{eq:vector-passive} and are given by
    \begin{gather} \label{eq:A-tilde}
      \tilde {\MAT A}_t = 
      \begin{bmatrix}
        \MAT A_t      & \MAT {B_tG_t}                        & \MAT {B_tH_t}                      \\ 
        \MAT {C_tA_t} & \MAT {C_tB_tG_t}                     & \MAT {C_tB_tH_t}                   \\ 
        \MAT 0        & \MAT P_{my,t}+\MAT P_{mu,t} \MAT G_t & \MAT P_{mm,t}+\MAT P_{mu,t} \MAT H_t
      \end{bmatrix},
      \\
      \label{eq:B-tilde}
      \tilde {\MAT B}_t =
       \begin{bmatrix}
         \MAT B_t \\ \MAT {C_tB_t} \\ \MAT P_{mu,t}
      \end{bmatrix},
      \quad\text{and}\quad
      \tilde {\MAT F}_t =
      \begin{bmatrix}
        \MAT I & \MAT 0 \\
        \MAT 0 & \MAT I \\
        \MAT 0 & \MAT 0
      \end{bmatrix};
    \end{gather}
    where the blocks in the first column of $\tilde {\MAT F}_t$ have dimensions
    compatible with $W^0_t$ and the blocks in the second column have dimensions
    compatible with $\VVEC(W^1_t, \dots, W^n_t)$. 
  \item 
    The observations are linear in the state and the control and may be written
    as
    \begin{align}
      \tilde Y_1 &= \MAT 0 \\
      \tilde Y_t &= \tilde {\MAT C}_t \tilde X_{t-1} + \tilde {\MAT D}_t \tilde
      U_{t-1},
      \quad t > 1
    \end{align}
    where $\tilde {\MAT C}_t$ and $\tilde {\MAT D}_t$ are matrices of
    appropriate dimensions given by
    \begin{gather}
      \tilde {\MAT C}_t = 
      \begin{bmatrix}
        \MAT 0 & \MAT P_{zy,t} + \MAT P_{zu,t} \MAT G_t &  \MAT P_{zm,t}+ \MAT P_{zu,t}\MAT H_t 
      \end{bmatrix},
        \\
        \tilde {\MAT D}_t = \MAT P_{zu,t}.
      \end{gather}
    \item 
      The per-step cost is quadratic in the state and control action and may be
      written as
      \begin{equation}
        \ell(X_t, U_t) = 
        \tilde \ell(\tilde X_t, \tilde U_t) = 
        \begin{bmatrix} 
          \tilde X_t^\TRANS & \tilde U_t^\TRANS 
        \end{bmatrix}
        \begin{bmatrix} 
          \tMAT Q_t & \tMAT N_t \\ \tMAT N_t^\TRANS & \tMAT R_t 
        \end{bmatrix}
        \begin{bmatrix}
          \tilde X_t \\ \tilde U_t 
        \end{bmatrix}
      \end{equation}
      where $\tMAT Q_t$, $\tMAT N_t$, $\tMAT R_t$ are obtained by
      combining~\eqref{eq:cost} and~\eqref{eq:vector-passive} and are given by
      \begin{gather}
        \tilde {\MAT Q}_t = 
        \begin{bmatrix}
          \MAT Q_t & \MAT 0 &\MAT 0 \\
          \MAT 0 & \MAT {G_t^\TRANS R_tG_t} &\MAT {G_t^\TRANS R_t H_t} \\
          \MAT 0 & \MAT {H_t^\TRANS R_t G_t} &  \MAT {H_t^\TRANS R_t H_t}
        \end{bmatrix}
        \\
        \tilde {\MAT N}_t =
        \begin{bmatrix}
          \MAT 0  \\
          \MAT {G_t^\TRANS R_t} \\
          \MAT {H_t^\TRANS R_t}
        \end{bmatrix}
        \quad\text{and}\quad
        \tilde {\MAT R}_t = \MAT R_t.
      \end{gather}
\end{enumerate}

Recall that we assume that $(\MAT G_{1:T}, \MAT H_{1:T})$ are fixed. The
auxiliary matrices $\tMAT A_t$, $\tMAT C_t$, $\tMAT Q_t$ and $\tMAT N_t$ defined
above depend on $\MAT G_t$ and $\MAT H_t$. 

\subsection{Characterization of the optimal controller}

The coordinated system defined above is a centralized partially observed LQG
system. Therefore, based on the standard results in linear stochastic
control~\cite{Caines:1987}, the optimal coordination strategy is characterized
as follows:

\begin{theorem} \label{thm:main-result}
  Define $\breve X_t$ as the estimate of the state $\tilde X_t$:
  \begin{align*}
    \breve X_t &= \EXP[ \tilde X_t \mid \tilde Y_{1:t}, \tilde U_{1:t-1} ]
  \end{align*}
  Then, we have
  \begin{enumerate}
    \item \emph{Kalman filtering update:} The initial value of the state
      estimate is given by 
      \(
        \breve X_1 = 0.
      \)
      For $t > 1$, the state estimate may be updated as follows
      \begin{align}
        \breve X_{t+1} = \tMAT A_t \breve X_t + \tMAT B_t \tilde U_t 
        &+ \tMAT A_t \tMAT P_t \tMAT C_t^\TRANS [ \tMAT C_t \tMAT P_t \tMAT C_t^\TRANS ]^{-1}
        (\tilde Y_{t+1} 
        - \tMAT C_{t+1} \breve X_t - \tMAT D_{t+1} \tilde U_t) \label{eq:coord_est}
      \end{align}
      where 
      \[
        \tMAT P_t = \EXP[ (\tilde X_t - \breve X_t)^2 \mid \tilde Y_{1:t}, \tilde U_{1:t-1} ],
      \]
      which may be computed a~priori by solving the following forward Riccati
      equation:
      \begin{align*}
        \tMAT P_1 &= \DIAG(\Sigma_x, \MAT 0_{d_y \times d_y}, \MAT 0_{d_m \times d_m})
        \\
        \tMAT P_{t+1} &= \tMAT A_t \tMAT P_t \tMAT A_t^\TRANS 
        + \tilde \Sigma_W \notag 
        - \tMAT A_t \tMAT P_t \tMAT C_t^\TRANS [\tMAT C_t \tMAT P_t \tMAT C_t^\TRANS]^{-1}
        \tMAT A_t^\TRANS \tMAT P_t \tMAT C_t
      \end{align*}
      where $d_y = \sum_{i=1}^n d_y^i$, $d_m = \sum_{i=1}^n d_m^i$, and
      $\tilde \Sigma_W$ is the covariance of $\tMAT F_t W_t$ which is
      given by
      \[
         \DIAG(\Sigma_{w^0}, \Sigma_{w^1}, \dots, \Sigma_{w^n}, \MAT 0)
      \]
      where $\MAT 0$ is a square matrix of dimension same as $M_t$. 

    \item \emph{Separation result:}
      The optimal action of the coordinator is given by
      \begin{equation}\label{eq:coord_action}
        \tilde U_t = \tMAT K_t \breve X_t
      \end{equation}
      where the gain matrices $\{\tMAT K_t\}_{t=1}^T$ are given by
      \begin{equation*}
        \tMAT K_t = - [\tMAT R_t + 
          \tMAT B_t^\TRANS \MAT S_{t+1} \tMAT B_t]^{-1} 
        \bm \Lambda_t
      \end{equation*}
      where
      \begin{equation*}
        \bm \Lambda_t = \tMAT N_t + \tMAT B_t^\TRANS \MAT S_{t+1} \tMAT A_t
      \end{equation*}
      and the matrices $\{\MAT S_t\}_{t=1}^T$ are given by backward Riccati
      equations:
      \begin{align*}
        \MAT S_T &= \tMAT Q_T \\
        \MAT S_t &= \tMAT A_t^\TRANS \MAT S_{t+1} \tMAT A_t + \tMAT Q_t \notag
        - \bm \Lambda_t^\TRANS 
        [\tMAT R_t + \tMAT B_t^\TRANS \MAT S_{t+1} \tMAT B_t]^{-1} 
        \bm \Lambda_t
      \end{align*}

    \item \emph{Performance:} The performance of the above strategy is given by
      \begin{equation*}
        J = 
        \sum_{\tau = 1}^{T} \TR[\tMAT P_t \tMAT Q_t + (\tilde \Sigma_W + \tMAT
          A_t \tMAT P_t \tMAT A_t^\TRANS - \tMAT P_{t+1}) \MAT S_{t+1}]
      \end{equation*}
  \end{enumerate}
\end{theorem}
Note that the matrices $(\tMAT K_{1:T}, \MAT S_{1:T}, \tMAT P_{1:T})$ obtained
above depend on the choice of the matrices $(\MAT G_{1:T}, \MAT H_{1:T})$. 

Since any linear strategy in the coordinated system  can be implemented in the  original system and vice versa, the above result gives the following structure of best linear strategies in the original system.
\begin{theorem}
\label{thm:main-result2}
  In Problem~(P1), the best linear control strategies are of the form
  \begin{align}
    U_t &= \tilde U_t +  \MAT G_t Y_t + \MAT H_t M_t \notag \\
    &= \tMAT K_t \breve X_t+ \MAT G_t Y_t + \MAT H_t M_t.
  \end{align}
  where $\MAT G_t = \DIAG(\MAT G^1_t, \ldots, \MAT G^n_t)$, $\MAT H_t = \DIAG(\MAT H^1_t, \ldots, \MAT H^n_t)$, 
  \begin{align*} 
  \breve X_t &= \VVEC(\hat X_t, \hat Y_t, \hat M_t) \\
  &=  \EXP[ \VVEC( X_t,  Y_t,  M_t)\mid Z_{1:t-1}, \tilde U_{1:t-1}],
 \end{align*}   
 and the evolution of $\breve X_t$,  the gain matrices $\tMAT K_t$ and the
 system performance are the same as in Theorem~\ref{thm:main-result}. 
\end{theorem}

\begin{remark}
  Let $\tMAT K_t = [\tMAT K^{1\,\TRANS}_t \mid \dots \mid \tMAT
  K^{n\,\TRANS}_t]^\TRANS_{\strut}$. Then, the control action of each
  controller may be written as
  \[
    U^i_t = \tMAT K^i_t \breve X_t + \MAT G^i_t Y^i_t + \MAT H^i_t M^i_t.
  \]
  Note that each controller is using its local information $(Y^i_t, M^i_t)$ and
  an estimate $\breve X_t$ based on the common information $Z_{1:t-1}$. 
\end{remark}

\begin{remark}
Note that for a given choice of $(\MAT G_{1:T}, \MAT H_{1:T})$, Theorem~\ref{thm:main-result} identifies the optimal $\tMAT K^i_t$ matrices and the associated cost. In order to find the best linear control strategies, we need to optimize the cost given in Theorem~\ref{thm:main-result} with respect to $(\MAT G_{1:T}, \MAT H_{1:T})$ --- which may be a non-convex optimization problem.
\end{remark}
\subsection{An equivalent representation of $\breve X_t$}\label{sec:equiv}

In Theorem~\ref{thm:main-result2}, it is possible to replace the estimate
$\breve X_t$ by a lower dimensional estimate $\breve S_t$ defined as
\[ 
  \breve S_t := \EXP[\VVEC( X_t,  M_t)\mid Z_{1:t-1}, \tilde U_{1:t-1} ].
\]

Given the definition of $\breve X_t$ in Theorem~\ref{thm:main-result}, we
immediately have that
\begin{equation}\label{eq:equiv_1}
\breve S_t = \begin{bmatrix} \MAT I & \MAT 0 & \MAT 0 \\ \MAT 0 & \MAT 0 & \MAT I\end{bmatrix} \breve X_t
\end{equation} 
Furthermore, since  $Y_t$ is related to $X_t$ through~\eqref{eq:observation}
and the primitive random variables are mutually independent, it follows that
\( \hat Y_t = \MAT C_t \hat X_t \) and 
therefore,
\begin{equation}\label{eq:equiv_2}
  \breve X_t = \begin{bmatrix}
   \MAT I   & \MAT 0 \\ 
   \MAT C_t & \MAT 0 \\ 
   \MAT 0   & \MAT I
  \end{bmatrix} \breve S_t
\end{equation} 

Equations~\eqref{eq:equiv_1} and \eqref{eq:equiv_2} imply that $\breve X_t$ can
be replaced by $\breve S_t$ as a sufficient statistic of common information in
Theorems~\ref{thm:main-result} and~\ref{thm:main-result2}. In particular, using
\eqref{eq:coord_action} and \eqref{eq:coord_est}, we get
\begin{equation}
 \tilde U_t = \tMAT K_t \begin{bmatrix}
   \MAT I & \MAT 0 \\
   \MAT C_t & \MAT 0\\
   \MAT 0 & \MAT I
\end{bmatrix} \breve S_t
\end{equation}
and
\begin{align}
        &\breve S_{t+1} = \begin{bmatrix} \MAT I & \MAT 0 & \MAT 0 \\ \MAT 0 & \MAT 0 & \MAT I\end{bmatrix}\Bigg[\tMAT A_t \begin{bmatrix}
   \MAT I & \MAT 0 \\
   \MAT C_t & \MAT 0\\
   \MAT 0 & \MAT I
\end{bmatrix} \breve S_t + \tMAT B_t \tilde U_t 
        +\notag \\ &\tMAT A_t \tMAT P_t \tMAT C_t^\TRANS [ \tMAT C_t \tMAT P_t \tMAT C_t ]^{-1}
        (\tilde Y_{t+1}  - \tMAT C_{t+1} \begin{bmatrix}
   \MAT I & \MAT 0 \\
   \MAT C_t & \MAT 0\\
   \MAT 0 & \MAT I
  \end{bmatrix} \breve S_t - \tMAT D_{t+1} \tilde U_t)\Bigg] \label{eq:S_update}
\end{align}

We can now state an equivalent version of Theorem \ref{thm:main-result2}.
\begin{theorem}
\label{thm:three}
  In Problem~(P1), the best linear control strategies are of the form
  \begin{align}
    U_t &= \tilde U_t +  \MAT G_t Y_t + \MAT H_t M_t \notag \\
    &= \tMAT L_t \breve S_t+ \MAT G_t Y_t + \MAT H_t M_t.
  \end{align}
  where $\MAT G_t = \DIAG(\MAT G^1_t, \ldots, \MAT G^n_t)$, $\MAT H_t = \DIAG(\MAT H^1_t, \ldots, \MAT H^n_t)$, 
  \begin{equation*} 
    \breve S_t = \VVEC(\hat X_t, \hat M_t) 
    =  \EXP[ \VVEC( X_t,   M_t)\mid Z_{1:t-1}, \tilde U_{1:t-1}],
  \end{equation*}   
  the evolution of $\breve S_t$ is given by \eqref{eq:S_update},  the gain matrices 
  \[
    \tMAT L_t = \tMAT K_t \begin{bmatrix}
      \MAT I   & \MAT 0 \\ 
      \MAT C_t & \MAT 0 \\ 
      \MAT 0   & \MAT I
    \end{bmatrix},
  \]
  the matrices $\tMAT K_t, \tMAT P_t$ and the system performance  are the same as in Theorem~\ref{thm:main-result}. 
\end{theorem}

\subsection{Generalization to models with common observations}

In some cases, in addition to the shared memory, controllers may also have a
common observation $Y^{\COM}_t$ about the state of the system given as
\begin{equation*}
  Y^{\COM}_t = \MAT C^{\COM}_t X_t + W^{\COM}_t,
\end{equation*}
where $W^{\COM}_t, t=1,2,\ldots,T$ is a sequence of i.i.d.\@ Gaussian variables
that are independent of the all the other primitive random variables. Each
controller can select its action according to a \emph{linear} control law of
the form
\begin{equation*}
  U^i_t = g^i_t(Y^i_t, M^i_t, \VEC Z_{1:t-1}, Y^{\COM}_{1:t}).
\end{equation*}
The methodology of Theorem~\ref{thm:main-result} can easily be adapted for
this model by allowing the coordinator to choose action $\tilde {\VEC U}_t =
\VVEC(\tilde U^1_t, \dots, \tilde U^n_t)$ based on the shared memory
\emph{and} the history of common observations. That is,
\begin{equation}
  \label{eq:coordinator-control2}
  \tilde U^i_t = \MAT K^i_t \VVEC(Z_{1:t-1}, Y^{\COM}_{1:t}).
\end{equation}
Following the same arguments as before, the coordinator's problem once again
becomes a classical LQG problem, thus establishing the result of
Theorem~\ref{thm:main-result} for this case with $\breve X_t$ now defined as
\begin{equation*}
  \breve X_t = \EXP[ (X_t, Y_t, M_t) \mid  Z_{1:t-1}, Y^{\COM}_{1:t}, \tilde U_{1:t-1}  ]
\end{equation*}

\subsection{Salient features of the result}

The above structural result shows that in the best linear strategy, the
control action at each time depends on the current local observation, the
current local memory, and a common information based estimate of the system
state and the local memories of all controllers. Thus, \important{the sufficient
  statistic is finite dimensional}.

Unlike prior work on structural results for decentralized control problems, our
result relies on the \emph{linearity} of the decentralized system and of control
strategies and not on partial nestedness or quadratic invariance. 

The result basically follows from two simple observations:  (i)~under linear
strategies, control actions can be viewed as superposition of two components---a local
information based component and a common information based component;  and
(ii)~once the matrices for calculating the local information based component have
been fixed, the problem of choosing the common information based component reduces to a centralized LQG problem.

\subsection{Comparison with \cite{NMT:partial-history-sharing}}

For decentralized control system with partial history sharing information
structure, it is shown in~\cite{NMT:partial-history-sharing} that the sufficient
statistic of the shared memory $Z_{1:t-1}$ is given by the posterior probability
distribution on $(X_t, M_t)$. In contrast, the result of Theorem~\ref{thm:three}
shows that when attention is restricted to linear strategies, the sufficient
statistic is given by the conditional mean $\breve S_t$ of $(X_t, M_t)$.
Therefore, the structural results of Theorem~\ref{thm:main-result} simplifies
the structural result of~\cite[Theorem~4]{NMT:partial-history-sharing} for LQG
systems with linear control strategies. 


Although the methodology used in proving Theorem~\ref{thm:main-result} and the
solution methodology of~\cite{NMT:partial-history-sharing} are similar, it is
not possible to derive the result of Theorem~\ref{thm:main-result} by directly
using the results of~\cite{NMT:partial-history-sharing}.
In~\cite{NMT:partial-history-sharing}, the coordinator solves a global
optimization problem to determine how controllers should use their local
information. On the other hand, to prove the result of
Theorem~\ref{thm:main-result}, we arbitrarily fix the components of the control
laws that use the local information   and then find the structure of the best
response strategies at the coordinator. 

This approach of fixing the part of control law that use the
local information and identifying the structure of coordinator's strategy was
also used for a two player partially nested problem in
\cite{LessardNayyar:2013}. In that paper, the authors used the  structure of
optimal linear strategies, along with the partial nestedness of the problem, to
explicitly derive the globally optimal control strategies. 

In contrast to the approach of~\cite{NMT:partial-history-sharing} which gives
the structure of globally optimal control laws and a dynamic programming
decomposition, our approach only gives the structure of best linear control
laws. It is not possible, in general, to extend our approach to find the best
linear control laws. The question whether the approach proposed in this paper
simplifies for partially nested teams warrants further investigation.

\section{Delayed Sharing Information Structure} \label{sec:examples}

In this section, we illustrate our results using  the specific example of
delayed sharing information structures. We consider two cases: (i) one with
symmetric delays where the observations and actions of any controller are
available to all other controllers after a delay of $k$ time steps and (ii) the
asymmetric delay case where the communication delay from controller $j$ to
controller $i$ is $k_{ij} < \infty$.

\subsection{Symmetric delays}

In delayed sharing information structure, each controller's observations and
control actions are shared with  all other controllers after a delay of $k \geq
1$ time steps~\cite{NMT:delay-sharing}. The system dynamics, local observations,
and cost function are the same as in Section~\ref{sec:model}. 

In the language of partial history sharing model, the shared memory in this case
consists of all observations and control actions that are at least $k $
time-steps old, that is, $\VEC Z_{1:t-1} = \VVEC(\VEC Y_{1:t-k}, \VEC
U_{1:t-k})$; and the local memory consists of the observations and actions taken at
$t-k+1$, \ldots, $t-1$, that is, $M^i_t = \VVEC(Y^i_{t-k+1:t-1},
U^i_{t-k+1:t-1})$. 

Therefore, the result of Theorem~\ref{thm:three} applies to this model with
\begin{itemize}
  \item the $\MAT P^i_{**}$ matrices in the memory update equations
    \eqref{eq:update_1} and \eqref{eq:update_2} are given by
    \begin{gather*}
      \allowdisplaybreaks
      \MAT P^i_{mm,t} = \begin{bmatrix}
        \MAT 0_{d^i_y \times d^i_m} \\
        \MAT 0_{d^i_u \times d^i_m} \\
        \MAT I_{(k-2)(d^i_y+d^i_u)} ~ \MAT 0_{(k-2)(d^i_y+d^i_u) \times (d^i_y +d^i_u)} 
      \end{bmatrix}
      \\
      \MAT P^i_{my,t} = \begin{bmatrix}
        \MAT I_{d^i_y} \\
        \MAT 0_{d^i_u \times d^i_y} \\
        \MAT 0_{(k-2)(d^i_y+d^i_u)\times d^i_y} 
      \end{bmatrix}
      \\
      \MAT P^i_{mu,t} = \begin{bmatrix}
        \MAT 0_{d^i_y \times d^i_u} \\
        \MAT I_{d^i_u} \\
        \MAT 0_{(k-2)(d^i_y+d^i_u)\times d^i_u}
      \end{bmatrix}
     \\
     \MAT P^i_{zm,t} = 
        \begin{bmatrix}
          \MAT 0_{(d^i_y + d^i_u) \times (k-2)(d^i_y+d^i_u)} ~ \MAT I_{(d^i_y + d^i_u)} 
        \end{bmatrix}
    \end{gather*}
    and $\MAT P^i_{zy,t} = \MAT 0$ and $\MAT P^i_{zu,t} = \MAT 0$.

  \item  and the estimate of Theorem \ref{thm:three} as 
    \[
      \breve S_t =  \EXP[ \VVEC( X_t,   M_t)\mid \VEC Y_{1:t-k}, \VEC U_{1:t-k}, \tilde U_{1:t-1}].
    \]
\end{itemize}

Recall that the evolution of the sufficient statistic $\breve S_t$ depends on
the choice of matrices  $(\MAT G_{1:T}, \MAT H_{1:T})$ in the control strategy. 
Such a dependence is also present in the sufficient statistic for
optimal control laws for the general delayed sharing
model~\cite{NMT:delay-sharing}. Hence, restricting attention to linear control
strategies does not lead to a two-way separation of estimation and control in
delayed sharing information structures. However, as we show next, it is possible
to have a one-way separation (estimation does not depend on control) if we keep
track of a subset of past observations and control actions at the coordinator.

\begin{corollary}
  \label{cor:delay-sharing-2}
  The result of Theorem \ref{thm:three} for the symmetric delay sharing model
  may be simplified as
  \begin{align}
    U_t &= \tilde U_t +  \MAT G_t Y_t + \MAT H_t M_t \notag \\
    &= \tMAT L_t  S_t+ \MAT G_t Y_t + \MAT H_t M_t.
  \end{align}
  where
  \begin{align*}
    &S_t = \VVEC( \hat X_{t-k+1|t-k}, \tilde{U}_{t-k+1:t-1}, Y_{t-2k+2:t-k}, U_{t-2k+2:t-k}), \notag \\
    &\hat X_{t-k+1|t-k} = \EXP[X_{t-k+1} \mid Y_{1:t-k},U_{1:t-k} ]
    \label{eq:cor2eq0} 
  \end{align*}
  and $\hat X_{t+1|t}$ is updated according to: 
  \begin{align*}
   \hat X_{1|0} &= 0 
   \\
   \hat X_{t+1|t} &= \MAT A_t \hat X_{t|t-1} + \MAT B_t  U_t 
     + \MAT A_t \MAT P_{t} \MAT C_t^\TRANS [\MAT C_t \MAT P_{t} \MAT C_t^\TRANS + \Sigma_w]^{-1}(Y_t-\MAT C_t \hat X_{t|t-1})
  \end{align*}
  where $\Sigma_w = \DIAG(\Sigma_{w^1}, \dots, \Sigma_{w^n})$ and $\MAT P_t = \EXP[
    (X_t - \hat X_{t|t-1})^2 \mid Y_{1:t-1}, U_{1:t-1}]$, which can be
    precomputed as follows:
  \begin{align*}
    \MAT P_1 &= \Sigma_x; \\
    \MAT P_{t+1} &= \MAT A_t \MAT P_t \MAT A_t^\TRANS +   \Sigma_{w^0} 
    - \MAT A_t \MAT P_{t} \MAT C_t^\TRANS [\MAT C_t \MAT P_{t} \MAT C_t^\TRANS + \Sigma_w]^{-1}\MAT C_t \MAT P_{t} \MAT A_t^\TRANS
  \end{align*}
\end{corollary}

See Appendix~\ref{app:cor} for a proof. Corollary~\ref{cor:delay-sharing-2}
shows that $S_t$ is a sufficient statistic for $(Y_{1:t-k}, U_{1:t-k})$. This
sufficient statistic consists of three parts:
\begin{enumerate}
  \item A strategy-independent $k$-step window $(Y_{t-2k+2:t-k},\allowbreak U_{t-2k+2:t-k})$
    of the history of observations and actions that are available to all controllers.
  \item A strategy-independent estimate of the $k$-step delayed state $X_{t-k+1}$
    based on the history of common information. Note that the update of $\hat
    X_{t-k+1|t-k}$ does not depend on the matrices $(\MAT G_{1:T}, \MAT H_{1:T})$. 
  \item A \emph{strategy-dependent} $k$-step window of the history of coordinated
    control actions $\tilde U_{t-k+1:t-1}$.
\end{enumerate}
This structure is similar to the optimal controller derived in~\cite[second
structural result]{NMT:delay-sharing}.

For the special case of delay $k=1$, the result of
Corollary~\ref{cor:delay-sharing-2} simplifies as follows.
\begin{corollary}
  \label{cor:delay-sharing-3}
  When the sharing delay $k=1$, the optimal control strategies may be chosen
  according to
  \begin{align}
    U_t &= \tilde U_t +  \MAT G_t Y_t  \notag \\
    &= \tMAT L_t  S_t+ \MAT G_t Y_t .
  \end{align}
  where
  \[
     \hat X_{t|t-1} = \EXP[X_{t} \mid Y_{1:t-1}, U_{1:t-1} ]
  \]
\end{corollary}
Corollary~\ref{cor:delay-sharing-3} is equivalent to the result obtained
in~\cite{Yoshikawa:1975, VaraiyaWalrand:1978}.

\subsection{Asymmetric delays}

In this model, controller $i$ observes the observations and control actions of
controller $j$ with a delay of $k_{ij} < \infty$. The information available to
controller $i$ at time $t$ consists of 
\[ 
  I^i_t = \{ Y^i_{1:t}, U^i_{1:t-1} \} \cup
  \bigcup_{j \neq i} \{Y^j_{1:t-k_{ij}}, U^j_{1:t-k_{ij}} \}.
\]

All delays are finite. For convenience, define $k_{ii}:=1$. Then, 
the information available to controller $i$ at time $t$ can be
written as
\[ 
  I^i_t = \{ Y^i_t \}  \cup
  \bigcup_{j = 1}^n \{Y^j_{1:t-k_{ij}}, U^j_{1:t-k_{ij}} \}.
\]

This information structure arises when controllers communicate along a
\emph{strongly connected graph with finite delay between any pair of
controllers}. The system dynamics, local observations, and cost function are the
same as in Section~\ref{sec:model}. Similar models have been considered
in~\cite{Rantzer:2007, LamperskiDoyle:2011, LamperskiDoyle:2012, FeyzmahdavianGattamiJohansson:2012}.
Note that unlike these models, we do not assume 
any sparsity structure on the matrices $\MAT A_t, \MAT B_t$ and $\MAT C_t$ in the system model.

\begin{figure}[!tbh]
  \centering
  \includegraphics[scale=0.8]{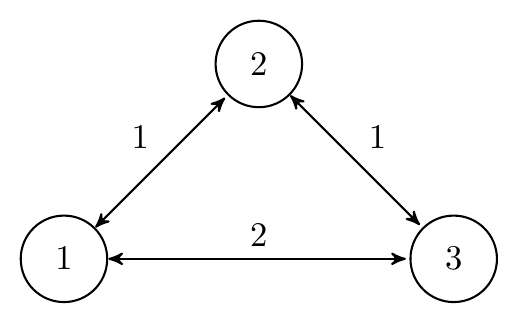}
  \caption{An example of a system with asymmetric delayed sharing. The number on
  the arrows denote the delay in flow of information.}
  \label{fig:example}
\end{figure}

Such a model has the generalized partial history sharing information structure.
As an illustration, consider the $3$ controller system shown in
Figure~\ref{fig:example}. Controllers~$1$ and~$2$ share information with
$1$-step delay, controllers~$2$ and~$3$ share information with $1$-step delay
but controllers~$1$ and~$3$ share information with $2$-step delay, that is,   
\begin{equation*}
  k_{12} = k_{21}= 1,\quad 
  k_{23} = k_{32}=1, \quad
  k_{13} = k_{31}=2.
\end{equation*}

The shared memory at time $t$ is given by
\[
  Z_{1:t-1} = \VVEC(Y^1_{1:t-2},U^1_{1:t-2},Y^2_{1:t-1},U^2_{1:t-1},Y^3_{1:t-2},U^3_{1:t-2});
\]
the local memories are
\begin{align*}
  M^1_t &= \VVEC(Y^1_{t-1},U^1_{t-1}), \\
  M^2_t &= \VVEC(Y^1_{t-1},U^1_{t-1},Y^3_{t-1},U^3_{t-1}),\\
  M^3_t &= \VVEC(Y^3_{t-1},U^3_{t-1});
\end{align*}
and the increment in shared memory at time $t$ is 
\[
  Z_t =\VVEC (Y^1_{t-1},U^1_{t-1},Y^2_t,U^2_t,Y^3_{t-1},U^3_{t-1}).
\]

The update of the local and shared memories may be written
as~\eqref{eq:update_3} and~\eqref{eq:update_4} with
\begin{gather*}
  \MAT P_{mm,t} = \MAT 0 \\
  \MAT P_{my,t} = \begin{bmatrix}
    \MAT I & \MAT 0 & \MAT 0 \\
    \MAT 0 & \MAT 0 & \MAT 0 \\
    \hline
    \MAT I & \MAT 0 & \MAT 0 \\
    \MAT 0 & \MAT 0 & \MAT 0 \\
    \MAT 0 & \MAT 0 & \MAT I \\
    \MAT 0 & \MAT 0 & \MAT 0 \\
    \hline
    \MAT 0 & \MAT 0 & \MAT I \\
    \MAT 0 & \MAT 0 & \MAT 0 \\
  \end{bmatrix},
  \quad
  \MAT P_{mu,t} = \begin{bmatrix}
    \MAT 0 & \MAT 0 & \MAT 0 \\
    \MAT I & \MAT 0 & \MAT 0 \\
    \hline
    \MAT 0 & \MAT 0 & \MAT 0 \\
    \MAT I & \MAT 0 & \MAT 0 \\
    \MAT 0 & \MAT 0 & \MAT 0 \\
    \MAT 0 & \MAT 0 & \MAT I \\
    \hline
    \MAT 0 & \MAT 0 & \MAT 0 \\
    \MAT 0 & \MAT 0 & \MAT I \\
  \end{bmatrix}
  \allowdisplaybreaks
  \\
  \MAT P_{zm,t} = \begin{bmatrix}
    \MAT I & \MAT 0 & \MAT 0 \\
    \hline
    \MAT 0 & \MAT 0 & \MAT 0 \\
    \hline
    \MAT 0 & \MAT 0 & \MAT I \\
  \end{bmatrix}
  \\
  \MAT P_{zy,t} = \begin{bmatrix}
    \MAT 0 & \MAT 0 & \MAT 0 \\
    \hline
    \MAT 0 & \MAT I & \MAT 0 \\
    \MAT 0 & \MAT 0 & \MAT 0 \\
    \hline
    \MAT 0 & \MAT 0 & \MAT 0 \\
  \end{bmatrix},
  \quad
  \MAT P_{zu,t} = \begin{bmatrix}
    \MAT 0 & \MAT 0 & \MAT 0 \\
    \hline
    \MAT 0 & \MAT 0 & \MAT 0 \\
    \MAT 0 & \MAT I & \MAT 0 \\
    \hline
    \MAT 0 & \MAT 0 & \MAT 0 \\
  \end{bmatrix}
\end{gather*}


Similar to the above example, the general model with asymmetric delays may be
considered as a special case of the generalized partial history sharing model.
For that matter, define $k^{*}_j := \max_{i} k_{ij}$. Thus, $k^{*}_j$ is the
delay after which controller $j$'s current information is available to all other
controllers. In the above example, $k^{*}_1=k^{*}_3=2$ and $k^{*}_2=1$.

Then, the common information available to all controllers at time $t$ is
\begin{equation*}
  Z_{1:t-1} = \VVEC(Y^1_{1:t-k^{*}_1},U^1_{1:t-k^{*}_1}\ldots,Y^n_{1:t-k^{*}_n},U^n_{1:t-k^{*}_n}),
\end{equation*}
and the local memory of controller $i$ is
\begin{align*}
  M^i_t &= I^i_t \setminus \{Y^i_t,Z_{1:t-1} \}\\
        &= \VVEC(Y^1_{t-k^{*}_1 +1:t-k_{i1}}, U^1_{t-k^{*}_1 +1:t-k_{i1}} ,\ldots, 
        Y^n_{t-k^{*}_n +1:t-k_{in}}, U^n_{t-k^{*}_n +1:t-k_{in}}).
\end{align*}

To facilitate writing the memory update equations of the form \eqref{eq:update_1} and
\eqref{eq:update_2} for the general asymmetric delay model, it is helpful to
define the following vectors:
\begin{equation}
  L^i_t = \VVEC(Y^i_{t-k_i^*+1:t-1}, U^i_{t-k_i^*+1:t-1}).
\end{equation}
$L^i_t$ denotes the observations and control actions of controller $i$ that have
not yet been shared with all controllers. $L^i_t$ takes values in $\reals^{d^i_l}$. $L^i_t$ is always a sub-vector of
$M^i_t$. Note that $L^i_t$ may be distinct from $M^i_t$ in general (see the
example above). More explicitly, the relation between $L^i_t$ and $M^i_t$ can be
written as 
\begin{gather}
  L^i_t = \begin{bmatrix}
    \MAT 0_{(k^*_i-1)(d^i_y+d^i_u) \times \sum_{j<i} (k^*_j-k_{ij})(d^j_y+d^j_u) } &
    \MAT I_{(k^*_i-1)(d^i_y+d^i_u)} &
    \MAT 0_{(k^*_i-1)(d^i_y+d^i_u) \times \sum_{j>i}  (k^*_j-k_{ij})(d^j_y+d^j_u) }
  \end{bmatrix} M^i_t \label{eq:LfromM}
\end{gather}


Define $L_t = \VVEC(L^1_t, \ldots,L^n_t)$. Note that $M^i_t$ is a sub-vector of
$L_t$. The explicit relation between $M^i_t$ and $L_t$ can be written as 
\begin{gather}
  M^i_t = \DIAG(J_{i1}, \ldots, J_{in}) L_t , 
 ~\mbox{where }
 J_{ij} = [\MAT 0_{(k^*_j-k_{ij})(d^j_y+d^j_u) \times (k_{ij}-1)(d^j_y+d^j_u)} ~ \MAT I_{(k^*_j-k_{ij})(d^j_y+d^j_u)}] \label{eq:MfromL}
\end{gather}

Furthermore, $L^i_t$ has an update equation similar to \eqref{eq:update_2}:
\begin{equation}
  L^i_{t+1}
    = \tMAT P^i_{\ell m} L^i_t +  \tMAT P^i_{\ell y}Y^i_t + \tMAT P^i_{\ell u} U^i_t    
  \label{eq:asymm_update}
\end{equation}
where
\begin{gather*}
  \tMAT P^i_{\ell m} =
  \begin{bmatrix}
    \MAT 0_{d^i_y \times d^i_l} \\
    \MAT 0_{d^i_u \times d^i_l} \\
    \MAT I_{(k^*_i-2)(d^i_y+d^i_u)} ~ \MAT 0_{(k^*_i-2)(d^i_y+d^i_u) \times (d^i_y +d^i_u)}
  \end{bmatrix}
  \\
  \tMAT P^i_{\ell y} =
    \begin{bmatrix}
      \MAT I_{d^i_y} \\
      \MAT 0_{d^i_u \times d^i_y} \\
      \MAT 0_{(k^*_i-2)(d^i_y+d^i_u)\times d^i_y}
    \end{bmatrix}
  \\
  \tMAT P^i_{\ell u} =
      \begin{bmatrix}
        \MAT 0_{d^i_y \times d^i_u} \\
        \MAT I_{d^i_u} \\
        \MAT 0_{(k^*_i-2)(d^i_y+d^i_u)\times d^i_u}
      \end{bmatrix}
\end{gather*}
The increment in shared memory can be written in terms of $L^i_t$ as
\begin{equation}
  Z^i_t
  =  \begin{bmatrix} \MAT 0_{(d^i_y + d^i_u) \times (k^*_i-2)(d^i_y+d^i_u)} ~ \MAT I_{(d^i_y + d^i_u)} \end{bmatrix}L^i_t  \label{eq:asymm_z}
\end{equation}

Therefore, the result of Theorem~\ref{thm:three} applies to this model with 
\begin{itemize}
\item The analogue of \eqref{eq:update_3} obtained by combining \eqref{eq:LfromM}, \eqref{eq:asymm_update} and \eqref{eq:MfromL}.
\item The analogue of \eqref{eq:update_4} obtained by combining \eqref{eq:LfromM} and \eqref{eq:asymm_z}.
\item and the estimate of Theorem \ref{thm:three} as 
\begin{align*}  \breve S_t 
    =  \EXP[ \VVEC( X_t,   M_t)\mid \VVEC(Y^1_{1:t-k^{*}_1},U^1_{1:t-k^{*}_1}\ldots,
      Y^n_{1:t-k^{*}_n},U^n_{1:t-k^{*}_n}), \tilde U_{1:t-1}].
\end{align*}
\end{itemize}

Analogous to Corollary~\ref{cor:delay-sharing-2}, we also have the following result in this model.
\begin{corollary}
  \label{cor:asym-delay-sharing-2}
  Define $k^* = \max_{i,j} k_{ij}$. The result of Theorem \ref{thm:three} for
  the asymmetric delay sharing model may be simplified as
  \begin{align}
    U_t &= \tilde U_t +  \MAT G_t Y_t + \MAT H_t M_t \notag \\
    &= \tMAT L_t  S_t+ \MAT G_t Y_t + \MAT H_t M_t.
  \end{align}
   where
  \begin{gather*}
    S_t  =     \VVEC( \hat X_{t-k^*+1|t-k^*}, \tilde{U}_{t-k^*+1:t-1}, Y_{t-2k^*+2:t-k^*},
    U_{t-2k^*+2:t-k^*}) 
    \notag\\
    \hat X_{t-k^*+1|t-k^*} = \EXP[X_{t-k^*+1} \mid Y_{1:t-k^*},U_{1:t-k^*} ]
  \end{gather*}
\end{corollary}
The proof is similar to the proof of Corollary \ref{cor:delay-sharing-2} in Appendix~\ref{app:cor}.

A $3$ controller system with asymmetric delays (in particular, $k_{21}=k_{32}=k_{13}=1$ and
$k_{12}=k_{23}=k_{31}=2$) and a partially nested information structure is considerd
in~\cite{FeyzmahdavianGattamiJohansson:2012}. The authors
of~\cite{FeyzmahdavianGattamiJohansson:2012} identify optimal control strategies
whose structural form is similar to our result above. Note that our results hold for \emph{any} strongly connected communication graph with finite delays.

\section{Models that reduce to partial history sharing} \label{sec:others}

The approach presented int his paper is also applicable to  models
that are not partial history sharing as such but can be reduced to one by using a
\emph{person-by-person approach}~\cite{MMRY:tutorial-CDC}. We illustrate this by means of
two examples presented below.
 
\subsection{Coupled subsystems with control sharing}

In the control sharing model considered in~\cite{M:control-sharing}\footnote
{The  model presented here is simpler than the model described in~\cite{M:control-sharing}.
The results also extend to the generalized models considered
in~\cite{M:control-sharing}, but we restrict attention to the more simpler model
for ease of exposition.}
the system consists on $n$-subsystems; each
subsystem has a co-located control station. Let $X^i_t$ denote the state of
subsystem~$i$ and $U^i_t$ the control action of controller~$i$. Let $X_t =
\VVEC(X^1_t, \dots, X^n_t)$ and $U_t = \VVEC(U^1_t, \dots, U^n_t)$. The system
dynamics are given by
\[
  X^i_{t+1} = \MAT A^i_t X^i_t + \MAT B^i_t U_t + W^i_t
\]
where $\MAT A^i_t$ and $\MAT B^i_t$ are matrices of appropriate dimensions. Note that the
next state of subsystem~$i$ depends on the current state of subsystem~$i$ and
the control actions of all controllers. The noise processes
$\{W^i_t\}_{t=1}^\infty$ are mutually independent and independent across time.
The cost is quadratic and given by~\eqref{eq:cost}.

Control station~$i$ observes the state of control station~$i$ and the one-step
delayed control actions of all controllers. Each controller has \emph{perfect
recall}. Therefore, action $U^i_t$ must be chosen based on the data $(X^i_{1:t},
U_{1:t-1})$. It is shown in~\cite[Proposition 3]{M:control-sharing} using a
person-by-person approach that there
is no loss of optimality in shedding $X^i_{1:t-1}$ and choosing $U^i_t$ based on
the data $(X^i_t, U_{1:t-1})$. We restrict attention to controllers that are linear
functions of this data, i.e., controllers for the form
\[
  U^i_t = \MAT K^i_t U_{1:t-1} + \MAT G^i_t X^i_t
\]

This model fits the general partial history sharing model described in
Section~\ref{sec:model} with
\begin{itemize}
  \item the local memory $M^i_t$ is empty;
  \item the local observation $Y^i_t$ is $X^i_t$;
  \item the shared memory $Z_{1:t-1}$ is  $U_{1:t-1}$
  \item the update of the shared memory given by~\eqref{eq:update-2}
    where $\MAT P^i_{zy}= \MAT 0$, $\MAT P^i_{zu} = \MAT I$ and $\MAT P_{**} =
    \DIAG(\MAT P_{**}^1, \dots, \MAT P_{**}^n)$.
\end{itemize}
The results of Theorem~\ref{thm:three} apply to this model with 
\[
  \breve S_t = \EXP[ X_t \mid U_{1:t-1} ].
\]

For this model, it is known that linear strategies are not globally optimal. The
optimal non-linear control strategy is given by the embedding of the
observations in the control actions~\cite{Bismut:1972}.

\subsection {One-sided one-step delayed sharing}

Consider two coupled subsystems with one-sided one-step delayed
sharing. Let $X^i_t$ denote the state of subsystem~$i$ and
$U^i_t$ denote the control action of subsystem~$i$. Let $X_t = \VVEC(X^1_t,
X^2_t)$ and $U_t = \VVEC(U^1_t, U^2_t)$. The dynamics are arbitrary and given
by~\eqref{eq:state}. At each time, controller~1 observes $\VVEC(X^1_t, X^2_{t-1})$: the
current state of subsystem~1 and the one-step delayed state of subsystem~2;
controller~2 observes $X^2_t$: the current state of subsystem~2. Thus,
controller~1 chooses its control actions based on the data $(X^1_{1:t},
U^1_{1:t-1}, X^2_{1:t-1}, U^2_{1:t-1})$ and controller~2 based on $(X^2_{1:t},
U^2_{1:t-1})$. 
The cost is quadratic and given by~\eqref{eq:cost}.

When $\MAT A$  and $\MAT B$ are lower block triangular, the model is partially
nested~\cite{HoChu:1972}. Such a model was considered
in~\cite{NayyarKalathilJain:2013}. A minor variation of this model (which was
also partially nested) was also considered in~\cite{SwigartLall:2010,KhorsandAlamGattami:2012}. 
The sparsity assumptions on $\MAT A$ and $\MAT B$ are needed to prove global
optimality of linear strategies; but, as we show below, not to identify the
sufficient statistics for linear strategies. 

The structure of controller~1 can be simplified by using a person-by-person
approach. For any arbitrary choice of control strategy
for controller~2, the subproblem of finding the \emph{best response} strategy at
controller~1 is a centralized stochastic control problem. It can be shown that
$(X^1_t, X^2_{1:t-1}, U^2_{1:t-1})$ is an information state of this subproblem.
Therefore, there is no loss of optimality in choosing $U^1_t$ based on the data
$(X^1_t, X^2_{1:t-1}, U^2_{1:t-1})$. We restrict attention to controllers that are linear
functions of the available data, i.e., controllers of the form
\begin{equation*}
  U^i_t = \MAT K^i_t \VVEC(X^2_{1:t-1}, U^2_{1:t-1}) + G^i_t X^i_t ; \\ 
\end{equation*}

This model fits the general partial history sharing model described in
Section~\ref{sec:model} with
\begin{itemize}
  \item the local memory $M^i_t$ is empty;
  \item the local observation $Y^i_t$ is $X^i_t$;
  \item the shared memory $Z_{1:t-1}$ is  $\VVEC(X^2_{1:t-1}, U^2_{1:t-1})$;
  \item the update of the shared memory given by~\eqref{eq:update-2}
    where 
    $\MAT P^1_{zy} = \MAT 0$, $\MAT P^2_{zu} = \MAT 0$, $\MAT P^2_{zy} = \MAT I$, 
    $\MAT P^2_{zu} = \MAT I$, and $\MAT P_{**} = \DIAG(\MAT P_{**}^1, \dots,
    \MAT P_{**}^n)$.
\end{itemize}
The results of Theorem~\ref{thm:three} apply to this model with 
\[
  \breve S_t = \EXP[ X_t \mid X^2_{1:t-1}, U^2_{1:t-1}, \tilde U_{1:t-1} ].
\]

The above structural result is similar to the result obtained
in~\cite{NayyarKalathilJain:2013}. However, unlike
\cite{NayyarKalathilJain:2013}, our model does not have a partially nested
information structure. This suggests that the structure of the best linear
control law is a consequence of the linearity of control strategies rather than the
partially nested information structure.

\section {Conclusion}

Linear control strategies for LQG systems are appealing due to their analytical
and implementation simplicity. However, to fully leverage the advantages of
linear strategies, we need to identify finite dimensional sufficient statistics
for best linear strategies that can be easily updated. We identified such a
result in Theorem~\ref{thm:three} for decentralized systems with partial
history sharing information structures. The result relied on the {linearity} of
the decentralized system and is applicable to models that are neither partially
nested nor quadratically invariant.
 
We focused on the partial history sharing model in this paper because it
provides a common model for decentralized systems where controllers' local
information remains finite dimensional but the common information 
increases with time. 

We showed that our results provide sufficient statistics for different
variations of delayed sharing information structures, including those with
asymmetric delays that arise when controllers communicate along a strongly
connected graph.

We also showed that our approach is applicable to some decentralized systems
where local information is also increasing with time, provided one can first employ a
person by person optimality approach to find a preliminary sufficient statistic
which ensures that local information is finite dimensional. 
 
We have focused only on finding the structure of best linear control strategies
in this paper. It is not possible, in general, to extend our approach to compute
the best linear control strategies. Even in the absence of a complete
methodology to find the best linear strategies, the structural results of
Theorem~\ref{thm:three} are useful because they restrict the solution
space to search for best linear strategies. Furthermore, as is the case with the
sufficient statistics in centralized stochastic control, the sufficient
statistics of Theorem~\ref{thm:three} allow us to formulate the problem of
finding and implementing the best linear control strategies over an infinite
horizon. 

\section*{Acknowledgments}

The authors thank Bension Kurtaran for suggesting this problem and
Demosthenis Teneketzis for helpful discussions.

\appendices

\section {Proof of Corollary~\ref{cor:delay-sharing-2}} \label{app:cor}

To prove the result, we will argue that $\breve S_t = \VVEC(\hat{X}_t, \hat{M}_t)$  is a linear function of $S_t$ for the symmetric delay sharing model. Therefore,
the control law of Theorem~\ref{thm:three} can be written in the form
specified in Corollary~\ref{cor:delay-sharing-2}.

Observe that according to the coordinated system dynamics in
\eqref{eq:coord_dynamics}, $(X_t, M_t)$ is a linear function of
$\tilde{X}_{t-k+1} = \VVEC(X_{t-k+1}, Y_{t-k+1}, M_{t-k+1})$,
$\tilde{U}_{t-k+1:t-1}$ and $W^0_{t-k+1:t-1}, W^{1:n}_{t-k+1:t-1}$. Therefore,
by linearity of conditional expectation, $(\hat{X}_t, \hat{M}_t)$ is a
linear function of the following three terms
\begin{enumerate}
  \item \mbox{$\EXP[ \VVEC(X_{t-k+1}, Y_{t-k+1}, M_{t-k+1}) \mid 
      Y_{1:t-k}, U_{1:t-k}, \tilde U_{1:t-k}]$.}

  \item $\EXP[ \tilde U_{t-k+1:t-1} \mid 
      Y_{1:t-k}, U_{1:t-k}, \tilde U_{1:t-k}]$.

  \item $\EXP[ W_{t-k+1:t-1} \mid 
      Y_{1:t-k}, U_{1:t-k}, \tilde U_{1:t-k}]$.
\end{enumerate}

Consider each of these terms separately. Recall that in delayed sharing
information structure $M_{t-k+1} = \VVEC(Y_{t-2k+2:t-k}, U_{t-2k+2:t-k})$ which
are included in the right hand side of conditioning in the first term.
Therefore, 
\begin{equation} \label{eq:app-1}
  \EXP[M_{t-k+1} \mid Y_{1:t-k},U_{1:t-k} ,\tilde {\VEC U}_{1:t-1}] 
  =
  \VVEC(Y_{t-2k+2:t-k}, U_{t-2k+2:t-k}).
\end{equation}
Furthermore, using~\eqref{eq:observation}
\begin{align}
  \EXP[Y^i_{t-k+1}\mid Y_{1:t-k},U_{1:t-k} ,\tilde {\VEC U}_{1:t-1}] 
  &=
  \MAT C^i \EXP[X_{t-k+1}\mid Y_{1:t-k},U_{1:t-k} ,\tilde {\VEC U}_{1:t-1}] 
  \notag \\
  &= \MAT C^i \EXP[X_{t-k+1}\mid Y_{1:t-k},U_{1:t-k}]
  \label{eq:app-2}
\end{align}
where we removed $\tilde{U}^i_{1:t-1}$ from the right hand side of conditioning
because it is a function of $(Y^i_{1:t-k},U^i_{1:t-k})$ which are included in the
right hand side of conditioning. 
Combining~\eqref{eq:app-1} and~\eqref{eq:app-2}, we get that
$\EXP[ \VVEC(X_{t-k+1}, Y_{t-k+1}, M_{t-k+1}) \mid Y_{1:t-k}, U_{1:t-k}, \tilde U_{1:t-k}]$
is a linear function of  $(\hat X_{t-k+1|t-k}, Y_{t-2k+2:t-k},
U_{t-2k+2:t-k})$, which is a sub-vector of $S_t$.

The second term $\EXP[ \tilde U_{t-k+1:t-1} \mid Y_{1:t-k}, U_{1:t-k}, \tilde
U_{1:t-k}]$ is simply $\tilde{U}_{t-k+1:t-1}$ which is also a sub-vector of
$S_t$. 

Since the primitive random variables are independent, the third term
$\EXP[ W_{t-k+1:t-1} \mid Y_{1:t-k}, U_{1:t-k}, \tilde U_{1:t-k}]$ is $0$.

Therefore, $\breve S_t = \VVEC(\hat{X}_t, \hat{M}_t)$  is a linear function of $S_t$, which implies the result of the corollary.

\bibliographystyle{IEEEtran}
\bibliography{IEEEabrv,../../collection,../../personal}

\begin{thebibliography}{10}
\providecommand{\url}[1]{#1}
\csname url@samestyle\endcsname
\providecommand{\newblock}{\relax}
\providecommand{\bibinfo}[2]{#2}
\providecommand{\BIBentrySTDinterwordspacing}{\spaceskip=0pt\relax}
\providecommand{\BIBentryALTinterwordstretchfactor}{4}
\providecommand{\BIBentryALTinterwordspacing}{\spaceskip=\fontdimen2\font plus
\BIBentryALTinterwordstretchfactor\fontdimen3\font minus
  \fontdimen4\font\relax}
\providecommand{\BIBforeignlanguage}[2]{{%
\expandafter\ifx\csname l@#1\endcsname\relax
\typeout{** WARNING: IEEEtran.bst: No hyphenation pattern has been}%
\typeout{** loaded for the language `#1'. Using the pattern for}%
\typeout{** the default language instead.}%
\else
\language=\csname l@#1\endcsname
\fi
#2}}
\providecommand{\BIBdecl}{\relax}
\BIBdecl

\bibitem{MMRY:tutorial-CDC}
A.~Mahajan, N.~Martins, M.~Rotkowitz, and S.~Y\"uksel, ``Information structures
  in optimal decentralized control,'' in \emph{Proc. 51st IEEE Conf. Decision
  and Control}, Maui, Hawaii, Dec. 2012, pp. 1291 -- 1306.

\bibitem{Witsenhausen:1968}
H.~S. Witsenhausen, ``A counterexample in stochastic optimum control,''
  \emph{{SIAM} Journal of Optimal Control}, vol.~6, no.~1, pp. 131--147, 1968.

\bibitem{LipsaMartins:2011b}
G.~M. Lipsa and N.~C. Martins, ``Optimal memoryless control in {Gaussian}
  noise: A simple counterexample,'' \emph{Automatica}, vol.~47, no.~3, pp.
  552--558, 2011.

\bibitem{Radner:1962}
R.~Radner, ``Team decision problems,'' \emph{Annals of Mathmatical Statistics},
  vol.~33, pp. 857--881, 1962.

\bibitem{HoChu:1972}
Y.-C. Ho and K.-C. Chu, ``Team decision theory and information structures in
  optimal control problems--{Part I},'' \emph{{IEEE} Trans. Autom. Control},
  vol.~17, no.~1, pp. 15--22, 1972.

\bibitem{Yuksel:2009}
S.~Y\"uksel, ``Stochastic nestedness and the belief sharing information
  pattern,'' \emph{{IEEE} Trans. Autom. Control}, pp. 2773--2786, Dec. 2009.

\bibitem{BamiehVoulgaris:2005}
B.~Bamieh and P.~Voulgaris, ``A convex characterization of distributed control
  problems in spatially invariant systems with communication constraints,''
  \emph{Systems and Control Letters}, vol.~54, no.~6, pp. 575--583, 2005.

\bibitem{RotkowitzLall:2006}
M.~Rotkowitz and S.~Lall, ``A characterization of convex problems in
  decentralized control,'' \emph{{IEEE} Trans. Autom. Control}, vol.~51, no.~2,
  pp. 274--286, 2006.

\bibitem{WhittleRudge:1974}
P.~Whittle and J.~Rudge, ``The optimal linear solution of a symmetric team
  control problem,'' \emph{Journal of Applied Probability}, vol.~11, no.~2, pp.
  377--381, Jun. 1974.

\bibitem{Yoshikawa:1975}
T.~Yoshikawa, ``Dynamic programming approach to decentralized stochastic
  control problems,'' \emph{{IEEE} Trans. Autom. Control}, vol.~20, no.~6, pp.
  796 -- 797, Dec. 1975.

\bibitem{VaraiyaWalrand:1978}
P.~Varaiya and J.~Walrand, ``On delayed sharing patterns,'' \emph{{IEEE} Trans.
  Autom. Control}, vol.~23, no.~3, pp. 443--445, 1978.

\bibitem{SwigartLall:2010}
J.~Swigart and S.~Lall, ``An explicit state-space solution for a decentralized
  two-player optimal linear-quadratic regulator,'' in \emph{American Control
  Conference (ACC), 2010}, 2010, pp. 6385--6390.

\bibitem{ShahParrilo:2013}
P.~Shah and P.~Parrilo, ``{${\cal H}_{2}$}-optimal decentralized control over
  posets: A state-space solution for state-feedback,'' \emph{{IEEE} Trans.
  Autom. Control}, vol.~58, no.~12, pp. 3084--3096, Dec. 2013.

\bibitem{LessardLall:2011}
L.~Lessard and S.~Lall, ``A state-space solution to the two-player
  decentralized optimal control problem,'' in \emph{Communication, Control, and
  Computing (Allerton), 2011 49th Annual Allerton Conference on}.\hskip 1em
  plus 0.5em minus 0.4em\relax IEEE, 2011, pp. 1559--1564.

\bibitem{KhorsandAlamGattami:2012}
O.~Khorsand, A.~Alam, and A.~Gattami, ``Optimal distributed controller
  synthesis for chain structures: Applications to vehicle formations,'' Apr.
  2012, arXiv:1204.1869.

\bibitem{NayyarKalathilJain:2013}
N.~Nayyar, D.~Kalathil, and R.~Jain, ``Optimal decentralized control with
  asymmetric one-step delayed information sharing,'' Sep. 2013,
  arXiv:1309.6376.

\bibitem{LessardNayyar:2013}
L.~Lessard and A.~Nayyar, ``Structural results and explicit solution for
  two-player {LQG} systems on a finite time horizon,'' \emph{arXiv preprint
  arXiv:1303.3256}, 2013.

\bibitem{NMT:partial-history-sharing}
A.~Nayyar, A.~Mahajan, and D.~Teneketzis, ``Decentralized stochastic control
  with partial history sharing: A common information approach,'' \emph{{IEEE}
  Trans. Autom. Control}, vol.~58, no.~7, pp. 1644--1658, Jul. 2013.

\bibitem{Nayyar:PhD}
A.~Nayyar, ``Sequential decision making in decentralized systems,'' Ph.D.
  dissertation, University of Michigan, Ann Arbor, MI, 2011.

\bibitem{NMT:delay-sharing}
A.~Nayyar, A.~Mahajan, and D.~Teneketzis, ``Optimal control strategies in
  delayed sharing information structures,'' \emph{{IEEE} Trans. Autom.
  Control}, vol.~56, no.~7, pp. 1606--1620, Jul. 2011.

\bibitem{Caines:1987}
P.~E. Caines, \emph{Linear stochastic systems}.\hskip 1em plus 0.5em minus
  0.4em\relax John Wiley \& Sons, Inc., 1987.

\bibitem{Rantzer:2007}
A.~Rantzer, ``A separation principle for distributed control,'' in
  \emph{Decision and Control, 2006 45th IEEE Conference on}, Dec. 2006, pp.
  3609--3613.

\bibitem{LamperskiDoyle:2011}
A.~Lamperski and J.~Doyle, ``On the structure of state-feedback lqg controllers
  for distributed systems with communication delays,'' in \emph{Decision and
  Control and European Control Conference (CDC-ECC), 2011 50th IEEE Conference
  on}, Dec. 2011, pp. 6901--6906.

\bibitem{LamperskiDoyle:2012}
------, ``Dynamic programming solutions for decentralized state-feedback lqg
  problems with communication delays,'' in \emph{American Control Conference
  (ACC)}, Jun. 2012, pp. 6322--6327.

\bibitem{FeyzmahdavianGattamiJohansson:2012}
H.~R. Feyzmahdavian, A.~Gattami, and M.~Johansson, ``Distributed
  output-feedback lqg control with delayed information sharing,'' in \emph{IFAC
  Workshop on Distributed Estimation and Control in Networked Systems
  (NECSYS)}, 2012.

\bibitem{M:control-sharing}
A.~Mahajan, ``Optimal decentralized control of coupled subsystems with control
  sharing,'' \emph{{IEEE} Trans. Autom. Control}, vol.~58, no.~9, pp.
  2377--2382, Sep. 2013.

\bibitem{Bismut:1972}
J.-M. Bismut, ``An example of interaction between information and control: The
  transparency of a game,'' \emph{{IEEE} Trans. Autom. Control}, vol.~18,
  no.~5, pp. 518--522, Oct. 1972.

\end{thebibliography}

\end{document}